\documentstyle{article}

\begin{document}
\baselineskip=18pt
\title{\hspace {11cm}{\small \bf IMPNWU-9709011}\\
\vspace {2cm} 
Dynamically twisted algebra $A_{q,p;\hat{\pi}}(\widehat{gl_2})$ 
as current algebra generalizing screening currents of q-deformed Virasoro 
algebra}

\author{
 Bo-yu Hou$^{b}$ and Wen-li Yang$^{a,b}$ 
\thanks{e-mail :wlyang@nwu.edu.cn}
\thanks{Fax    :0086-029-8303511}
\\
\bigskip\\
$^{a}$ CCAST ( World Laboratory ), P.O.Box 8730 ,Beijing 100080, China\\
$^{b}$ Institute of Modern Physics, Northwest 
University, Xian 710069, China
\thanks{Mailing address}}
\maketitle

\begin{abstract}
In this paper, we propose an elliptic algebra 
$A_{q,p;\hat{\pi}}(\widehat{gl_2})$  
 which is based on the relations $RLL=LLR^{*}$, where $R$ and $R^{*}$ 
 are the dynamical R-maxtrices of $A^{(1)}_{1}$ type face model with the  
 elliptic moduli shifted by the center of the algebra.
 From the Ding-Frenkel correspondence , we find that its corresponding
 (Drinfeld) current algebra at level
 one is the algebra of screening currents for q-deformed Virasoro algebra.
We realize the elliptic algebra at level one by Miki's construction from
the bosonization for the type I and type II vertex operators.
We also show that the algebra
$A_{q,p;\hat{\pi}}(\widehat{gl_2})$ is related with the algebra
$A_{q,p}(\widehat{gl_2})$ by a dynamically twisting.

\end{abstract}
\section{Introduction}
As the quantum form of fundamental Poisson bracket, the $``RLL=LLR"$
relations (or $``RLL"$ formalism)
define various quantum algebra which appear in quantum field theory (QFT) and
statistical mechanics.It is associated with  structure constants---R-matrix
satisfying the Yang-Baxter equation.Drinfeld and Jimbo have discovered a
fundamental algebra structure--- quantum algebra $U_q(g)$ where $g$ is some
finite or infinite dimensional Lie algebra[3,4].
Faddeev , Reshetikhin and Takhtajan[5]  realize the algebra $U_q(g)$ 
(FRT construction),where g is some finite dimensional Lie algebra,  
by the ``RLL" formalism with spectra parameter independent R-matrix.
Later, Reshetikhin  and Semenov-Tian-Shansky[6]
 constructed a new realization of q-deformed affine algebra by the ``RLL" 
formalism with trigonometric R-matrix (which is so-called RS relations)
characterized by the spectra parameter shifted with the center of the algebra.
Ding and Frenkel[7] gave the isomorphism between the realization given by 
Reshetikhin  and Semenov-Tian-Shansky and the Drinfeld realization of 
q-deformed affine algebra. Moreover, Khoroshkin[8] constructed successfully
the realization of Yangian Double with center $DY_{\hbar}(\widehat{g})$ in
the``RLL"  formalism[9], which is associated with rational R-matrix.
Foda et al proposed an elliptic extension of the quantum affine 
algebra $A_{q,p}(\widehat{sl_2})$[10] as an symmetric algebra for the 
eight-vertex model.The elliptic algebra is based on   
the generalized ``RLL" formalism : $RLL=LLR^{*}$ , where $R$ and $R^{*}$ 
are eight-vertex R-matrices with elliptic moduli differing by an amount 
depending on the level k of the representation on which $L$ acts. The 
algebra $A_{\hbar ,\eta}(\widehat{sl_2})$ as the scaling limit of the  
elliptic algebra $A_{q,p}(\widehat{sl_2})$ was formulated in $``RLL=LLR^{*}"$
formulation by Jimbo et al[36] and was studied by Khoroshikin
et al through the Gauss decomposition[8].In fact, the  above progress in
``RLL" formalism which we mentioned are all involved  in vertex type models.

Another progress of quantum algebra focus on the q-deformation [11-13]
and $\hbar$-deformation[14] of 
Virasoro and W algebra, and  q-deformed extented Virasoro algebra[27],
which would play quite the same role in
off-critical integrable model as that of Virasoro ,W algebra and extented
Virasoro algebra[38-40] 
in two-dimensional conformal field theories (CFT)[15] for the 
critical model. Q-deformed Virasoro (q-Virasoro) algebra and q-deformed
extented Virasoro(extented q-Virasoro) algebra arise in
the two-dimensional solvable lattice models [12,17,27]  (e.g ABF
model[19] etc.);   
$\hbar$-deformed Virasoro ($\hbar$-Virasoro) algebra is studied as 
the hidden symmetry of the massive integrable field models 
(e.g. the Restricted sine-Gordon model [14]).It was shown that the screening  
currents for q-Virasoro algebra[12,13,20,21] and $\hbar$-Virasoro algebra
[14] satisfy a closed algebra relations which is some further deformation of
q-affine algebra  and Yangian double with center. In another way , 
q-Virasoro algebra, extented q-Virasoro and $\hbar$-Virasoro algebra
can be redefined as the
algebra which commute with the screening currents up to a total 
difference[12-14,17,27]. Apparently, they constitute the hidden symmetries
in $A^{(1)}_{1}$ type face  model[12] .
 So, the studies of the algebraic structure of the screening 
currents for q-Virasoro, extented q-Virasoro and
$\hbar$-Virasoro  are of great importance.In
this paper, we mainly  deal with the screening algebra for q-Virasoro 
algbera and possibly the extented q-Virasoro ,as a byproduct,
the screening algebra for $\hbar$-Virasoro algebra and related
$\hbar$-deformed extented Virasoro algebra
can be obtained by taking scaling limit of the q-deformed version.

The ``RLL" formalism was originally formulated for the nondynamical 
Yang-Baxter equation (or the vertex type models), Felder[25]
succeeded in extending it to incorporate dynamical Yang-Baxter equations 
(Gervais-Neveu-Felder equation)[22] which associte with q-deformation of 
Knizhnik-Zamolodchikov-Bernard equation 
on torus. In fact, the ``RLL" formalism given by  Felder[25] and
Enriquez et al[22] is a dynamical version of the FRT construction and
RS relations respectively.
In this paper, we extend the works of Fodal et al[10] to the dynamical
R-matrix case.
Namely, we proposed an elliptic algebra $A_{q,p;\hat{\pi}}(\widehat{gl_2})$ 
based on the relations $RLL=LLR^{*}$ , where $R$ and $R^{*}$ are the 
dynamical R-matrices for $A_{1}^{(1)}$ face model (i.e the solution to 
Star-Triangle relation in $A_{1}^{(1)}$ type face model) with elliptic 
moduli shifted by the center of the algebra. Using the Ding-Frenkel
correspondence, we construct  Drinfeld currents for the algebra 
$A_{q,p;\hat{\pi}}(\widehat{gl_2})$.From the Drinfeld currents for the 
algebra (which is a subalgebra of $A_{q,p;\hat{\pi}}(\widehat{gl_2})$),
we show that the (Drinfeld) current algebra at level one and a higer level is 
just the algebra of screening currents for q-Virasoro algebra  and
the algebra of screening currents for the extented q-Virasoro algebra[27] 
respectively.
The algebra of screening currents at level one 
was studied by Awata et al[12] for q-Virasoro algebra and by Feigin et 
al[13,21]  for q-deformed W algebra, which is some elliptic deformation  
of affine algebra. The elliptic algebra 
$A_{q,p;\hat{\pi}}(\widehat{gl_2})$ at higher level would 
play an important role in the studies of the fusion ABF 
models[26,29] and relate with the extended q-Virasoro algebra.
Moreover, the algebra  $A_{q,p;\hat{\pi}}(\widehat{gl_2})$ is the 
dynamical twisted algebra[23-25] of the elliptic algebra 
$A_{q,p}(\widehat{gl_2})$.

The paper is organized as follows.In section 2, after some review of  
q-Virasoro algebra, we introduce the algebra of screening currents 
for q-Virasoro algbera . In section 3, we construct  an
elliptic algbera $A_{q,p;\hat{\pi}}(\widehat{gl_2})$ in terms of
$L^{\pm}$-operator
which satisfy the dynamical relations of $RLL=LLR^{*}$ formulation.From 
 Ding-Frenkel correspondence, we find that the corresponding dynamical
 Drinfeld current form a subalgebra which structure constants do not
 depend upon the dynamical variable and   is just the algebra
 of screening currents  for extented q-Virasoro algebra.
 The twisted relations between the algebra
 $A_{q,p;\hat{\pi}}(\widehat{gl_2})$ and $A_{q,p}(\widehat{gl_2})$ 
 is constructed. In section 4, the bosonization of the type I
 and type II vertex operator for the algebra $A_{q,p;\hat{\pi}}(\widehat{gl_2})$ 
 at level one is constructed. By the Miki's construction, we obtained 
 the bosonization for the algebra $A_{q,p;\hat{\pi}}(\widehat{gl_2})$ at level one.
 The corresponding generalizing algebra
 $A_{\hbar,\eta;\hat{\pi}}(\widehat{gl_2})$ for $\hbar$-deformed Virasoro
 algebra and $\hbar$-deformed extented Virasoro algebra ,which is the scaling
 limit of the algebra $A_{q,p;\hat{\pi}}(\widehat{gl_2})$, is studied in
 section 5. Finally, we
 give summary and discussions in section 6. Appendix contains some 
 detailed calculations.

\section{Algebra of screening currents for q-Virasoro algebra}
We start with defining q-Virasoro algebra and corresponding quantum 
Miura transformation
\subsection{q-Virasoro algebra and quantum Miura transformation}
Let $w$ be a generic complex number with $Im(w)>0$ and $r$ be a real number 
with $4<r$, and set $x=e^{i\pi w}$. Define elliptic functions 
\begin{eqnarray*}
& &\theta \left[\begin{array}{c}a\\b\end{array}\right](z,\tau)
=\sum_{m\in Z}exp\{i\pi[(m+a)^2\tau +2(m+a)(z+b)]\}\ \ ,\ \ Im(\tau)>0\\
& &\sigma_{\alpha}=\sigma_{(\alpha_1,\alpha_2)}(z,\tau)=
\theta \left[\begin{array}{c}\frac{1}{2}+\frac{\alpha_1}{2}
\\\frac{1}{2}+\frac{\alpha_2}{2}\end{array}\right](z,\tau)\\
& &\theta^{(k)}(z,\tau)=\theta \left[\begin{array}{c}-\frac{k}{2}\\0
\end{array}\right](z,2\tau)
\end{eqnarray*}
We will use the following abbreviation
\begin{eqnarray*}
& &[v]_{t}=x^{\frac{v^{2}}{t}-v}\Theta_{x^{2t}}(x^{2v})\ \ ,\ \ 
=\sigma_{0}(\frac{v}{t},-\frac{1}{tw})\times const.\ \ ,\ \ 0<t\\
& &\Theta_{q}(z)=(z;q)(qz^{-1};q)(q;q)\ \ \ ,\ \ \ 
(z;q_1,...,q_m)=\prod_{i_1,...i_m=0}^{\infty}(1-zq^{i_1}...q^{i_m})
\end{eqnarray*}
\noindent Q-Virasoro algebra generated by $\{T(z)\}$ with the following 
relations[11-13]
\begin{eqnarray}
f(\frac{w}{z})T(z)T(w)-f(\frac{z}{w})T(w)T(z)=\frac{(x^{r}-x^{-r})
(x^{(r-1)}-x^{-(r-1)})}{x-x^{-1}}(\delta(\frac{w}{x^2z})-
\delta(\frac{x^2w}{z}))
\end{eqnarray}
\noindent where $\delta(z)=\sum_{n\in Z}z^{m}$ and 
\begin{eqnarray*}
f(z)=(1-z)^{-1}\frac{(zx^{2r};x^4)(zx^{2-2r};x^4)}
{(zx^{2+2r};x^4)(zx^{4-2r};x^4)}
\end{eqnarray*}
\noindent The generators $T(z)$ for q-Virasoro algebra can be obtained 
by the following quantum Miura transformation
\begin{eqnarray}
T(z)=\Lambda(x^{-1}z)+\Lambda^{-1}(xz)
\end{eqnarray}
\noindent Define q-deformed bosonic oscillators 
$\beta_{m}\ \ (m\in Z/\{0\})$ 
\begin{eqnarray}
[\beta_{m},\beta_{n}]=m\frac{(x^m-x^{-m})(x^{(r-1)m}-x^{-(r-1)m})}
{(x^{2m}-x^{-2m})(x^{rm}-x^{-rm})}\delta_{m+n,0}
\end{eqnarray}
\noindent and zero mode operator $P$ and $Q$ such that $[P,iQ]=1$

Then the fundamental operator $\Lambda(z)$ can be realized by  q-deformed 
bosonic oscillators Eq.(3) as follows
\begin{eqnarray}
\Lambda(z)=x^{\sqrt{2r(r-1)}P}:exp\{-\sum_{m\neq}^{\infty}(x^{rm}-x^{-rm})
\frac{\beta_{m}}{m}z^{-m}\}:
\end{eqnarray}

\subsection{Algebra of screening currents}
Let us introduce the screening currents  $E(v),F(v)$ for 
q-Virasoro algebra
\begin{eqnarray}
& &E(v)=e^{i\sqrt{\frac{2(r-1)}{r}}(Q-i2vlnxP)}:e^{\sum_{m\neq 0}
\frac{x^{m}+x^{-m}}{m}\beta_{m}x^{-2vm}}:\\
& &F(v)=e^{-i\sqrt{\frac{2r}{r-1}}(Q-i2vlnxP)}:e^{-\sum_{m\neq 0}
\frac{x^{m}+x^{-m}}{m}\beta'_{m}x^{-2vm}}:
\end{eqnarray}
\noindent where $\beta'_{m}=\frac{x^{rm}-x^{-rm}}{x^{(r-1)m}-x^{-(r-1)m}}
\beta_{m}$.Besides the well-known bosonic realization of screening currents
$E(v),F(v)$, let us introduce $H^{\pm}$
\begin{eqnarray}
& &H^{-}(v)=-x^{-4v}:E(v+\frac{1}{4})F(v-\frac{1}{4}):\\
& &H^{+}(v)=-x^{-4v}:E(v-\frac{1}{4})F(v+\frac{1}{4}):
\end{eqnarray}
The screening currents commute with the generators of q-Virasoro
algebra  up to a total difference and form a closed algebra. In fact,
from the normal order in appendix A., one can find that the screening currents defined 
in Eq.(5)---Eq.(8) realize an algebra satisfying the relations
\begin{eqnarray}
& &E(v_1)E(v_2)=\frac{[v_1-v_2-1]_{r}}{[v_1-v_2+1]_{r}}E(v_2)E(v_1)\\
& &F(v_1)F(v_2)=\frac{[v_1-v_2+1]_{r-1}}{[v_1-v_2-1]_{r-1}}F(v_2)F(v_1)\\
& &[E(v_1),F(v_2)]=\frac{1}{x-x^{-1}}\{
\delta(v_1-v_2+\frac{1}{2})H^{+}(v_1+\frac{1}{4})
-\delta(v_1-v_2-\frac{1}{2})H^{-}(v_1-\frac{1}{4})\}\\
& &H^{\pm}(v_1)E(v_2)=\frac{[v_1-v_2-1\mp\frac{1}{4}]_{r}}
{[v_1-v_2+1\mp\frac{1}{4}]_{r}}E(v_2)H^{\pm}(v_1)\\
& &H^{\pm}(v_1)F(v_2)=\frac{[v_1-v_2+1\pm\frac{1}{4}]_{r-1}}
{[v_1-v_2-1\pm\frac{1}{4}]_{r-1}}F(v_2)H^{\pm}(v_1)\\
& &H^{\pm}(v_1)H^{\pm}(v_2)=\frac{[v_1-v_2+1]_{r-1}[v_1-v_2-1]_{r}}
{[v_1-v_2-1]_{r-1}[v_1-v_2+1]_{r}}H^{\pm}(v_2)H^{\pm}(v_1)\\
& &H^{+}(v_1)H^{-}(v_2)=\frac{[v_1-v_2+1+\frac{1}{2}]_{r-1}
[v_1-v_2-1-\frac{1}{2}]_{r}}
{[v_1-v_2-1+\frac{1}{2}]_{r-1}[v_1-v_2+1-\frac{1}{2}]_{r}}
H^{-}(v_2)H^{+}(v_1)
\end{eqnarray}
\noindent where $H^{-}(v)=H^{+}(v+\frac{1}{2}-r)$.The algebra of screening
currents written by Awata[12] has the similar algebraic  relations to ours.
Actually, the screening currents
defined in Eq.(9)---Eq.(15) realize the (Drinfeld) current algebra of
an elliptic algebra $A_{q,p;\hat{\pi}}(\widehat{gl_{2}})$ at level one, 
which will be given by ``RLL" formalism in the following section .

\section{The dynamical algebra $A_{q,p;\hat{\pi}}(\widehat{gl_{2}})$}
We propose an elliptic algebra $A_{q,p;\hat{\pi}}(\widehat{gl_2})$ based on a 
dynamical $RLL=LLR^{*}$ relations.Then by Ding-Frenkel correspondence, we
will show that its Drinfeld current algebra is related to 
the current algebra generalizing the screening currents for
q-Virasoro algebra. Moreover, algebra  $A_{q,p;\hat{\pi}}(\widehat{gl_2})$ 
is an algebraic structure underlying the elliptic solution to the 
Star-Triangle relation in $A^{(1)}_{1}$face type model including ABF[19] and its
fused version[26,27,29].

\subsection{The R-matrix}
Define a dynamical elliptic R-matrix ( R-matrix for $A^{(1)}_{1}$ face  
model[28,30])
\begin{eqnarray}
& &R_{F}(v,\hat{\pi})\equiv R_{F}(v,\hat{\pi},r)=\left(
\begin{array}{llll}a&&&\\&b&c&\\&d&e&\\&&&a\end{array}
\right)
\end{eqnarray}
\noindent where $\hat{\pi}$ is the dynamical variable corresponding the 
height for the face type model and enjoy in some relations with  algebra  
$A_{q,p;\hat{\pi}}(\widehat{gl_{2}})$ (See Eq.(27) and Eq.(28)). 
The matrix elements of  the R-matrix is defined by 
\begin{eqnarray}
& &a(v,\hat{\pi})=x^{\frac{1-r}{r}v}\frac{g_{1}(v)}{g_{1}(-v)}\ \ ,\ \ 
g_{1}(v)=\frac{\{x^{2+2v}\}\{x^{2+2r+2v}\}}{\{x^{4+2v}\}\{x^{2r+2v}\}}
\ \ ,\ \ {z}=(z;x^{2r},x^{4})\\
& &\frac{b(v,\hat{\pi})}{a(v,\hat{\pi})}=\frac
{[v]_{r}[\hat{\pi}-1]_{r}}{[v+1]_{r}[\hat{\pi}]_{r}}\ \ ,\ \ 
\frac{c(v,\hat{\pi})}{a(v,\hat{\pi})}=\frac
{[v+\hat{\pi}]_{r}[1]_{r}}{[v+1]_{r}[\hat{\pi}]_{r}}\\
& &\frac{d(v,\hat{\pi})}{a(v,\hat{\pi})}=\frac
{[\hat{\pi}-v]_{r}[1]_{r}}{[v+1]_{r}[\hat{\pi}]_{r}}\ \ ,\ \ 
\frac{e(v,\hat{\pi})}{a(v,\hat{\pi})}=\frac
{[v]_{r}[1+\hat{\pi}]_{r}}{[v+1]_{r}[\hat{\pi}]_{r}}
\end{eqnarray}
\noindent One can see that $a(v,\hat{\pi})$ does not depend on the 
dynamical variable $\hat{\pi}$. Moreover, let us introduce two R-matrices 
$R_{F}^{\pm}$ which coincide with $R_{F}$ in Eq.(16) up to scalar factors 
independent on the dynamical variable
\begin{eqnarray}
& &R_{F}^{\pm}(v,\hat{\pi})\equiv R_{F}^{\pm}(v,\hat{\pi},r)
=\tau^{\pm}(v)R_{F}(v,\hat{\pi})\ \ ,\ \ 
\tau^{\pm}(v)=\tau(-v\pm \frac{1}{2})\\
& &\tau(v)=x^{-v}\frac{(x^{1+2v};x^{4})(x^{3-2v};x^{4})}
{(x^{3+2v};x^{4})(x^{1-2v};x^{4})}\nonumber
\end{eqnarray}
\noindent where $\tau^{\pm}(v)$ are the same as that of Foda et al [10].
$R_{F}^{\pm}$ are regarded as linear operators on $V\otimes V$, 
with $V=span\{e_{\pm}\}$. Let h be the diagonal $2\times 2$ matrix  
Diag(1,-1) such that $he_{\pm}=\pm e_{\pm}$. The dynamical R-matrices
$R_{F}^{\pm}(v,\hat{\pi})$ satisfy
the dynamical Yang-Batxer equation (i.e the  modified Yang-Baxter 
equation [22]) in $V\otimes V\otimes V$
\begin{eqnarray}
& &R^{\pm}_{F12}(v_1-v_2,\hat{\pi}-2h^{(3)})R^{\pm}_{F13}
(v_1-v_3,\hat{\pi})R^{\pm}_{F23}(v_2-v_3,\hat{\pi}-2h^{(1)})\nonumber\\
& &\ \ \ \ =R^{\pm}_{F23}(v_2-v_3,\hat{\pi})
R^{\pm}_{F13}(v_1-v_3,\hat{\pi}-2h^{(2)})R^{\pm}_{F12}
(v_1-v_2,\hat{\pi})
\end{eqnarray}
\noindent Here we choose the same notation as Enriquez et al[22]: 
$R^{\pm}_{F12}(v,\hat{\pi}-2h^{(3)})$ means that if 
$a\otimes b\otimes e_{\mu} \in V\otimes V\otimes V ,\mu\in \pm$, then 
$R^{\pm}_{F12}(v,\hat{\pi}-2h^{(3)})a\otimes b\otimes e_{\mu}$ = 
$R^{\pm}_{F12}(v,\hat{\pi}-2\mu)a\otimes b\otimes e_{\mu}$, and
the other symbols have a similar meaning. Besides the dynamical Yang-Baxter equation Eq.(21), 
the R-matrices have the following properties:
\begin{eqnarray}
& &{\rm \large \bf Unitarity}\ \ \ \ \ R^{\pm}_{F12}(v,\hat{\pi})
R^{\mp}_{F21}(-v,\hat{\pi})=id\ \ \ \ \ \ \ \ \ \ \\
& & {\rm \large \bf Crossing \ \  relations}\ \ \ 
R_{F}^{\pm}(-1-v,\hat{\pi})^{\mu'\nu'}_{\mu\nu}=\mu\mu'
R_{F}^{\mp}(v,\hat{\pi}-\mu')^{-\nu\mu'}_{-\nu'\mu}\frac{[\hat{\pi}-\mu']_{r}}
{[\hat{\pi}]_{r}}
\end{eqnarray}
\noindent Moreover, the R-matrices $R_{F}^{\pm}(v,\hat{\pi})$ enjoys in the
following property
$$
R^{+}_{F}(v+r,\hat{\pi})=R^{-}_{F}(v,\hat{\pi})\eqno(23a)
$$
Here, if $R^{\pm}_{F12}(v,\hat{\pi})=\sum a_{i}\otimes b_{i}$,
with $a_{i},b_{i}\in {\rm End}(V)$, then  
$R^{\pm}_{F21}(v,\hat{\pi})=\sum b_{i}\otimes a_{i}$.
\subsection{The algebra $A_{q,p;\hat{\pi}}(\widehat{gl_2})$}
Let us proceed to the definition of the elliptic algebra 
$A_{q,p;\hat{\pi}}(\widehat{gl_2})$. Consider $L^{\pm}$-operators
\begin{eqnarray*}
L^{\pm}(v,\hat{\pi})=\left(\begin{array}{ll}
L^{\pm}_{++}&L^{\pm}_{+-}\\L^{\pm}_{-+}&L^{\pm}_{--}\end{array}\right)
\end{eqnarray*}
\noindent which matrix elements are the generators of the  elliptic algebra 
$A_{q,p;\hat{\pi}}(\widehat{gl_2})$ given by the following 
commutation relations:
\begin{eqnarray}
& &R_{F}^{+}(v_1-v_2+\frac{c}{2},\hat{\pi})L^{+}_{1}(v_1,\hat{\pi})
L^{-}_{2}(v_2,\hat{\pi})=L^{-}_{2}(v_2,\hat{\pi})
L^{+}_{1}(v_1,\hat{\pi})R_{F}^{*+}(v_1-v_2-\frac{c}{2},\hat{\pi})\\
& &R_{F}^{-}(v_1-v_2-\frac{c}{2},\hat{\pi})L^{-}_{1}(v_1,\hat{\pi})
L^{+}_{2}(v_2,\hat{\pi})=L^{+}_{2}(v_2,\hat{\pi})
L^{-}_{1}(v_1,\hat{\pi})R_{F}^{*-}(v_1-v_2+\frac{c}{2},\hat{\pi})\\
& &R_{F}^{\pm}(v_1-v_2,\hat{\pi})L^{\pm}_{1}(v_1,\hat{\pi})
L^{\pm}_{2}(v_2,\hat{\pi})=L^{\pm}_{2}(v_2,\hat{\pi})
L^{\pm}_{1}(v_1,\hat{\pi})R_{F}^{*\pm}(v_1-v_2,\hat{\pi})
\end{eqnarray}
\noindent where $L^{\pm}_{1}(v,\hat{\pi})=L^{\pm}(v,\hat{\pi})\otimes id$,
 $L^{\pm}_{2}(v,\hat{\pi})=id\otimes L^{\pm}(v,\hat{\pi})$, 
 $R_{F}^{*\pm}(v,\hat{\pi})=R_{F}^{\pm}(v,-\hat{\pi},r-c)$ and c is the 
 center of the algebra (its vaule on some representation of the algebra 
 is usual called it as level of the algebra) . Moreover, 
 the $L^{\pm}$-operators are related with the dynamical variable $\hat{\pi}$: 
\begin{eqnarray*}
\hat{\pi} L^{(\pm)\mu}_{\ \ \ \nu}(v,\hat{\pi})=
L^{(\pm)\mu}_{\ \ \ \nu}(v,\hat{\pi})(\hat{\pi}+(\nu r-(r-c)\mu))
\end{eqnarray*}
\noindent Hence , the dynamical R-matrices have the following properties
\begin{eqnarray}
& &R_{F}^{\epsilon}(v,\hat{\pi})L^{(\epsilon')\mu}_{\ \ \ \nu}(v,\hat{\pi})
=L^{(\epsilon')\mu}_{\ \ \ \nu}(v,\hat{\pi})R_{F}^{\epsilon}
(v,\hat{\pi}+\mu c)\\
& &R_{F}^{*\epsilon}(v,\hat{\pi})L^{(\epsilon')\mu}_{\ \ \ \nu}(v,\hat{\pi})
=L^{(\epsilon')\mu}_{\ \ \ \nu}(v,\hat{\pi})R_{F}^{*\epsilon}(v,\hat{\pi}+\nu c)
\end{eqnarray}
\noindent where $\epsilon ,\epsilon' \in \pm$ and the following property 
\begin{eqnarray*}
[v+t]_{t}=-[v]
\end{eqnarray*}
\noindent is used.

\noindent {\bf Remark:} The relation Eq.(25) is the direct result of 
Eq.(22) and Eq.(24) .

Let
\begin{eqnarray*}
L^{\pm}(v,\hat{\pi})=
\left(\begin{array}{cc}1&0\\E^{\pm}(v)&1\end{array}\right)
\left(\begin{array}{cc}K^{\pm}_{1}(v)&0\\&K^{\pm}_{2}(v)\end{array}\right)
\left(\begin{array}{cc}1&F^{\pm}(v)\\0&1\end{array}\right)
\end{eqnarray*}
\noindent be the Gauss decomposition of $L^{\pm}$-operators. 
For the convenience, we introduce the following symbols 
\begin{eqnarray*}
& &R_{F}^{\pm}(v,\hat{\pi})=\left(\begin{array}{llll}
a^{\pm}(v)&&&\\&b^{\pm}(v)&c^{\pm}(v)&\\&d^{\pm}(v)&e^{\pm}(v)&
\\&&&a^{\pm}(v)\end{array}
\right)\ \ ,\ \  R_{F}^{*\pm}(v,\hat{\pi})=\left(\begin{array}{llll}
a'^{\pm}(v)&&&\\&b'^{\pm}(v)&c'^{\pm}(v)&\\&d'^{\pm}(v)&e'^{\pm}(v)&
\\&&&a'^{\pm}(v)\end{array}
\right)  
\end{eqnarray*}
\noindent {\bf Remark:} The elements $a^{\pm}(v)$ and $a'^{\pm}(v)$ 
do not depend on the dynamical variable , and commute with the Gauss  
components of $L^{\pm}$-operators. 

\noindent Define the total currents $E(v)$ and  
$F(v)$ by the corresponding Ding-Frenkel correspondence
\begin{eqnarray}
& &E(v)=E^{+}(v)-E^{-}(v+\frac{c}{2})\ \ \ ,\ \ \ 
F(v)=F^{+}(v+\frac{c}{2})-F^{-}(v)
\end{eqnarray}
\noindent Then we have the following

\noindent {\bf Proposition 1.} The total currents $E(v),F(v)$ and 
$K^{\pm}_{i}(v)$ $(i=1,2)$ satisfy the following commutation relations 
$$
a^{\pm}(v_1-v_2)K^{\pm}_{i}(v_1)K^{\pm}_{i}(v_2)=
K^{\pm}_{i}(v_2)K^{\pm}_{i}(v_1)a'^{\pm}(v_1-v_2) \eqno(30a)$$
$$
b^{\pm}(v_1-v_2)K^{\pm}_{1}(v_1)K^{\pm}_{2}(v_2)=
K^{\pm}_{2}(v_2)K^{\pm}_{1}(v_1)b'^{\pm}(v_1-v_2) \eqno(30b)$$
$$
a^{+}(v_1-v_2+\frac{c}{2})K^{+}_{i}(v_1)K^{-}_{i}(v_2)=
K^{-}_{i}(v_2)K^{+}_{i}(v_1)a'^{\pm}(v_1-v_2-\frac{c}{2}) \eqno(30c)$$
$$
b^{+}(v_1-v_2+\frac{c}{2})K^{+}_{1}(v_1)K^{-}_{2}(v_2)=
K^{-}_{2}(v_2)K^{+}_{1}(v_1)b'^{+}(v_1-v_2-\frac{c}{2}) \eqno(30d)$$
$$
b^{-}(v_1-v_2-\frac{c}{2})K^{-}_{1}(v_1)K^{+}_{2}(v_2)=
K^{+}_{2}(v_2)K^{-}_{1}(v_1)b'^{-}(v_1-v_2+\frac{c}{2}) \eqno(30e)$$
$$
K^{+}_{1}(v_1)E(v_2)K^{+}_{1}(v_1)^{-1}=\frac{a^{+}(v_1-v_2)}
{b^{+}(v_1-v_2)}E(v_2)\eqno(31a)$$
$$
K^{+}_{2}(v_2)E(v_1)K^{+}_{2}(v_2)^{-1}=E(v_1)\frac{a^{+}(v_1-v_2)}
{b^{+}(v_1-v_2)}\eqno(31b)$$
$$
K^{-}_{1}(v_1)E(v_2)K^{-}_{1}(v_1)^{-1}=\frac{a^{+}(v_1-v_2-\frac{c}{2})}
{b^{+}(v_1-v_2-\frac{c}{2})}E(v_2)\eqno(31c)$$
$$
K^{-}_{2}(v_2)E(v_1)K^{-}_{2}(v_2)^{-1}=E(v_1)\frac{a^{+}(v_1-v_2+\frac{c}{2})}
{b^{+}(v_1-v_2+\frac{c}{2})}\eqno(31d)$$
$$
K^{+}_{1}(v_1)^{-1}F(v_2)K^{+}_{1}(v_1)=F(v_2)\frac{a'^{+}(v_1-v_2-\frac{c}{2})}
{b'^{+}(v_1-v_2-\frac{c}{2})}\eqno(32a)$$
$$
K^{+}_{2}(v_2)^{-1}F(v_1)K^{+}_{2}(v_2)=\frac{a'^{+}(v_1-v_2+\frac{c}{2})}
{b'^{+}(v_1-v_2+\frac{c}{2})}F(v_1)\eqno(32b)$$
$$
K^{-}_{1}(v_1)^{-1}F(v_2)K^{-}_{1}(v_1)=F(v_2)\frac{a'^{+}(v_1-v_2)}
{b'^{+}(v_1-v_2)}\eqno(32c)$$
$$
K^{-}_{2}(v_2)^{-1}F(v_1)K^{-}_{2}(v_2)=\frac{a'^{+}(v_1-v_2)}
{b'^{+}(v_1-v_2)}F(v_1)\eqno(32d)$$
$$
E(v_1)\frac{a^{\pm}(v_1-v_2)}{b^{\pm}(v_1-v_2)}E(v_2)
=E(v_2)\frac{a^{\pm}(v_2-v_1)}{b^{\pm}(v_2-v_1)}E(v_1)\eqno(33a)$$
$$
F(v_1)\frac{a'^{\pm}(v_2-v_1)}{b'^{\pm}(v_2-v_1)}F(v_2)
=F(v_2)\frac{a'^{\pm}(v_1-v_2)}{b'^{\pm}(v_1-v_2)}F(v_1)\eqno(33b)$$
\addtocounter{equation}{4}
\begin{eqnarray}
& &[E(v_1),F(v_2)]=\frac{1}{x-x^{-1}}\{\delta (v_2-v_1-\frac{c}{2}) 
K^{-}_{2}(v_1+\frac{c}{2})\frac{[\hat{\pi}]_{r-c}[1]_{r-c}}
{\theta'_{r-c}[\hat{\pi}-1]_{r-c}}K^{-}_{1}(v_1+\frac{c}{2})^{-1}
\nonumber\\
& &\ \ \ \ \ \ \delta (v_2-v_1+\frac{c}{2})
K^{+}_{2}(v_1)\frac{[\hat{\pi}]_{r-c}[1]_{r-c}}
{\theta'_{r-c}[\hat{\pi}-1]_{r-c}}K^{+}_{1}(v_1)^{-1}\}
\end{eqnarray}
\noindent where 
$K^{-}_{i}(v)=K^{+}_{i}(v-\frac{c}{2}+r)$ and 
$\theta'_{t}=(x-x^{-1})\frac{\partial}{\partial v}[v]_{t}|_{v=0}$.

\noindent The proof of these relations is shifted to the appendix B.

\noindent {\bf Remark:} The elliptic algebra $U_{q,p}(\widehat{sl_2})$ 
proposed by Konno[27], has  different $K(v)$ from ours($K^{\pm}_{i}(v)$
also in Ref.[37]),
in which the commutation relations between $K(v)$ and $E(v),F(v)$ 
do not depend on the dynamical variable.However, they share the same
subalgebra generated by $H^{\pm}(v)$,$E(v)$ and $F(v)$ (see below) .

Set
\begin{eqnarray}
& &H^{+}(v)
=K^{-}_{2}(v+\frac{c}{4})\frac{[\hat{\pi}]_{r-c}[1]_{r-c}}
{\theta'_{r-c}[\hat{\pi}-1]_{r-c}}K^{-}_{1}(v+\frac{c}{4})^{-1}\\
& &H^{-}(v)
=K^{+}_{2}(v+\frac{c}{4})\frac{[\hat{\pi}]_{r-c}[1]_{r-c}}
{\theta'_{r-c}[\hat{\pi}-1]_{r-c}}K^{+}_{1}(v+\frac{c}{4})^{-1}
\end{eqnarray}
\noindent Like the case of q-affine  algebra[7] and 
Yangian double algebra[8,9], we can obtain the corresponding Drinfeld
current algebra of $A_{q,p;\hat{\pi}}(\widehat{gl_2})$ which is the
subalgebra of the elliptic algebra
$A_{q,p;\hat{\pi}}(\widehat{gl_2})$ and generated by $E(v),F(v),
H^{\pm}(v)$. Namely, we have 

\noindent {\bf Proposition 2.} The ellipitc (Drinfeld) current algebra of 
algebra $A_{q,p;\hat{\pi}}(\widehat{gl_2})$ is generated
by $E(v),F(v),H^{\pm}(v)$ with the following algebraic relations
\begin{eqnarray}
& &E(v_1)E(v_2)=\frac{[v_1-v_2-1]_{r}}{[v_1-v_2+1]_{r}}E(v_2)E(v_1)\\
& &F(v_1)F(v_2)=\frac{[v_1-v_2+1]_{r-c}}{[v_1-v_2-1]_{r-c}}F(v_2)F(v_1)\\
& &[E(v_1),F(v_2)]=\frac{1}{x-x^{-1}}\{
\delta(v_1-v_2+\frac{c}{2})H^{+}(v_1+\frac{c}{4})
-\delta(v_1-v_2-\frac{c}{2})H^{-}(v_1-\frac{c}{4})\}\\
& &H^{\pm}(v_1)E(v_2)=\frac{[v_1-v_2-1\mp\frac{c}{4}]_{r}}
{[v_1-v_2+1\mp\frac{c}{4}]_{r}}E(v_2)H^{\pm}(v_1)\\
& &H^{\pm}(v_1)F(v_2)=\frac{[v_1-v_2+1\pm\frac{c}{4}]_{r-c}}
{[v_1-v_2-1\pm\frac{c}{4}]_{r-c}}F(v_2)H^{\pm}(v_1)\\
& &H^{\pm}(v_1)H^{\pm}(v_2)=\frac{[v_1-v_2+1]_{r-c}[v_1-v_2-1]_{r}}
{[v_1-v_2-1]_{r-c}[v_1-v_2+1]_{r}}H^{\pm}(v_2)H^{\pm}(v_1)\\
& &H^{+}(v_1)H^{-}(v_2)=\frac{[v_1-v_2+1+\frac{c}{2}]_{r-c}
[v_1-v_2-1-\frac{c}{2}]_{r}}
{[v_1-v_2-1+\frac{c}{2}]_{r-c}[v_1-v_2+1-\frac{c}{2}]_{r}}
H^{-}(v_2)H^{+}(v_1)
\end{eqnarray}
\noindent and
\begin{eqnarray*}
H^{-}(v)=H^{+}(v+\frac{c}{2}-r)
\end{eqnarray*}
\noindent The proof of these formulas is shifted to the appenix C.

\noindent {\bf Remark:} 1.The deformed parameters are $q=x$ and $p=x^{2r}$
(cf. Ref.[10]).Moreover, the constructure coefficient in Eq.(37)---Eq.(43) do
not depend on the dynamical variable $\hat{\pi}$.

$\ \ \ \ \ \ \ \ \ $2.When $r\longrightarrow +\infty$, the limit current
algebra is the algebra $U_{q}(\widehat{sl_2})$[32].

One can see that if c=1 (i.e level one ), the current algebra of algebra  
$A_{q,p;\hat{\pi}}(\widehat{gl_2})$ be the algebra of screening currents 
for q-Virasoro algebra (cf. Eq.(9)---Eq.(15)) 
which play the role of symmetry algebra in ABF model[12,17].
For the general level k $(k\in$ integer),the 
algebra $A_{q,p;\hat{\pi}}(\widehat{gl_2})$ would correpond to the k-fusion 
ABF model[26,27,29] and in this case, some q-deformation of the extended 
Virasoro algebra[38-40] would exist in such a way that their screening 
currents satisfy current algebra of algebra
$A_{q,p;\hat{\pi}}(\widehat{gl_2})$.
So, this elliptic algebra would play an important role in the studies 
of $A^{(1)}_{1}$ type face models as that of the algebra  
$A_{q,p}(\widehat{gl_2})$ in the eight-vertex model[10].

\subsection{The algebra $A_{q,p;\hat{\pi}}(\widehat{gl_{2}})$ as the 
dynamically twisted algebra $A_{q,p}(\widehat{gl_{2}})$}
It is well-known that there exists a face-vertex correspondence between 
$A^{(1)}_1$ face model and eight-vertex model when r is a generic one
[26,28,29].This would result in the $``$equivalence" between the underlying 
algebra $A_{q,p;\hat{\pi}}(\widehat{gl_{2}})$ and 
$A_{q,p}(\widehat{sl_{2}})$
--- the algebra $A_{q,p;\hat{\pi}}(\widehat{gl_{2}})$ is the dynamcially  
twisted algebra of $A_{q,p}(\widehat{gl_{2}})$. In this section, we
restrict our attention to the case of r being a generic one.

Let $\epsilon_{\mu}$ $(\mu \in \pm)$ be the orthonormal basis in $R^2$, 
which are supplied with the inner product $<\epsilon_{\mu},\epsilon_{\nu}>
=\delta_{\mu\nu}$. Set
\begin{eqnarray*}
\overline{\epsilon}_{\mu}=\epsilon_{\mu}-\epsilon\ \ ,\ \ \epsilon=\frac{
\epsilon_{-}+\epsilon_{+}}{2}
\end{eqnarray*}
\noindent Then , define the intertwiners [28,33]
\begin{eqnarray*}
& &\varphi^{(m)}_{\hat{k},\mu}(v)=\theta^{(m)}(\frac{v+
<\hat{k},\overline{\epsilon}_{\mu}>}{r},
-\frac{1}{rw})\\
& &\varphi'^{(m)}_{\mu,\hat{l}}(v)=\theta^{(m)}(\frac{v+
<\hat{l},\overline{\epsilon}_{\mu}>}{r-c},
-\frac{1}{(r-c)w})\\
& &<\hat{k},\overline{\epsilon}_{\mu}>\equiv \mu\hat{k}\ \ \ ,\ \ \ 
<\hat{l},\overline{\epsilon}_{\mu}>\equiv \mu\hat{l} 
\end{eqnarray*}
$$
\hat{\pi}\equiv (r-c)\hat{k}-r\hat{l}\eqno(43a)
$$
\noindent Here, we remark that the decomposition of Eq.(43a) can be defined
 only for the generic r[33] and
\begin{eqnarray*}
\hat{k}L^{(\pm)\mu}{\nu}(v,\hat{\pi})=
L^{(\pm)\mu}{\nu}(v,\hat{\pi})(\hat{k}+\mu c)\ \ ,\ \ 
\hat{l}L^{(\pm)\mu}{\nu}(v,\hat{\pi})=
L^{(\pm)\mu}{\nu}(v,\hat{\pi})(\hat{l}+\nu c)
\end{eqnarray*}
\noindent The face-vertex correspondence relations read as
\begin{eqnarray*}
& &R^{ij}_{mn}(v_1-v_2)\varphi^{(m)}_{\hat{k},\nu}(v_1)
\varphi^{(n)}_{\hat{k}-\overline{\epsilon}_{\nu},\mu}(v_2)
=\sum_{\mu'\nu'}R^{\ \ \nu\mu}_{F\nu'\mu'}(v_1-v_2,\hat{\pi})
\varphi^{(i)}_{\hat{k}-\overline{\epsilon}_{\mu'},\nu'}(v_1)
\varphi^{(j)}_{\hat{k},\mu'}(v_2)\\
& &R^{*ij}_{mn}(v_1-v_2)\varphi'^{(m)}_{\nu,\hat{l}+
\overline{\epsilon}_{\mu}}(v_1)\varphi'^{(n)}_{\mu,\hat{l}}(v_2)
=\sum_{\mu\nu}R^{*\nu\mu}_{F\nu'\mu'}(v_1-v_2,\hat{\pi})\varphi'^{(i)}_{\nu',
\hat{l}}(v_1)
\varphi'^{(j)}_{\mu',\hat{l}+\overline{\epsilon}_{\nu'}}(v_2)
\end{eqnarray*}
\noindent where the nondynamical R-matrices $R$ and $R^{*}$ are the same 
as that of Foda et al [10].
Moreover, we can introduce intertwiners $\overline{\varphi}_
{\hat{k},\mu}$ and $\overline{\varphi}'_{\mu,\hat{l}}$ satisfying 
relations[28]
\begin{eqnarray*}
& &\sum_{m}\overline{\varphi}^{(m)}_{\mu,\hat{k}}
\varphi^{(m)}_{\nu,\hat{k}}=\delta_{\mu,\nu}\ \ ,\ \ \ 
\sum_{\mu}\overline{\varphi}^{(i)}_{\mu,\hat{k}}
\varphi^{(j)}_{\mu,\hat{k}}=\delta^{ij}\\ 
& &\sum_{m}\overline{\varphi}'^{(m)}_{\hat{l},\mu}
\varphi'^{(m)}_{\hat{l},\nu}=\delta_{\mu,\nu}\ \ ,\ \ \ 
\sum_{\mu}\overline{\varphi}'^{(i)}_{\hat{l},\mu}
\varphi'^{(j)}_{\hat{l},\mu}=\delta^{ij} 
\end{eqnarray*}
\noindent Then, we have the twisted relations between the R-matrix of 
eight-vertex model and the R-matrix of $A^{(1)}_{1}$ face model
\begin{eqnarray}
& &R^{\ \ \nu\mu}_{F\nu'\mu'}(v_1-v_2,\hat{\pi})=\sum_{ijmn}
\overline{\varphi}^{(j)}_{\hat{k},\mu'}(v_2)
\overline{\varphi}^{(i)}_{\hat{k}-\overline{\epsilon}_{\mu'},\nu'}(v_1)
R^{ij}_{mn}(v_1-v_2)\varphi^{(m)}_{\hat{k},\nu}(v_1)
\varphi^{(n)}_{\hat{k}-\overline{\epsilon}_{\nu},\mu}(v_2)\\
& &R^{*\nu\mu}_{F\nu'\mu'}(v_1-v_2,\hat{\pi})=\sum_{ijmn}
\overline{\varphi}'^{(j)}_{\mu',\hat{l}+\overline{\epsilon}_{\nu'}}(v_2)
\overline{\varphi}'^{(i)}_{\nu',\hat{l}}(v_1)R^{*ij}_{mn}(v_1-v_2)
\varphi'^{(m)}_{\nu,\hat{l}+\overline{\epsilon}_{\mu}}(v_1)
\varphi'^{(n)}_{\mu,\hat{l}}(v_2)
\end{eqnarray}
\noindent Moreover, we can constructed the twisted relations between 
the corresponding $L^{\pm}$-operators
\begin{eqnarray}
& &L^{\pm\nu}_{\ \ \mu}(v,\hat{\pi})=\sum_{mm'}
\overline{\varphi}^{(m')}_{\hat{k},\nu}(v)
L^{\pm m'}_{\ \ m}(v)\varphi'^{(m)}_{\mu,\hat{l}}(v)\nonumber\\
& &L^{\pm m'}_{\ \ m}(v)=\sum_{\mu\nu}\varphi^{(m')}_{\hat{k},\nu}(v)
L^{\pm\nu}_{\ \ \mu}(v,\hat{\pi})\overline{\varphi}'^{(m)}_{\mu,\hat{l}}(v)
\end{eqnarray}
\noindent Then we have 

\noindent {\bf Proposition 3.} The $L^{\pm}(v)$ operators given by the 
twisted relations Eq.(46) satisfying the commutation relations of  
the algebra $A_{q,p}(\widehat{gl_2})$[10] 
\begin{eqnarray}
& &R^{+}(v_1-v_2+\frac{c}{2})L^{+}_{1}(v_1)L^{-}_{2}(v_2)=
L^{-}_{2}(v_2)L^{+}_{1}(v_1)R^{*+}(v_1-v_2-\frac{c}{2})\\
& &R^{\pm}(v_1-v_2)L^{\pm}_{1}(v_1)L^{\pm}_{2}(v_2)=
L^{\pm}_{2}(v_2)L^{\pm}_{1}(v_1)R^{*\pm}(v_1-v_2)
\end{eqnarray}
\noindent where the r-matrices $R^{\pm}(v),R^{*\pm}(v)$ are the same as
that of Foda et al
\begin{eqnarray*}
R^{\pm}(v)\equiv \tau^{\pm}(v)R(v)\ \ \ ,\ \ \ R^{*\pm}(v)=
R^{*\pm}(v)|_{r\longrightarrow r-c}
\end{eqnarray*}

\section{The type I and Type II vertex and Miki's construction}
This section is devoted to the realization of an infinite dimensional 
representations of the algebra $A_{q,p;\hat{\pi}}(\widehat{gl_2})$  at 
level one by  the q-primary fields of q-Virasoro algebra.
\subsection{The type I and type II vertex operators}
The method of bosonization provides a powerful method to study the solvable 
lattice model both in vertex type model[32] and the face type model[17,30,35]. 
In this subsection , we give the bosonization of the type I[17] and typ II[35] 
vertex operator in ABF model by one free field. 

The type I vertex operator corresponds to the half-column transfer matrix 
 of the model, and type II vertex operator is expected to creat the 
 eigenstates of the transfer matrix. We denote the two types of vertex 
 operator as 
\begin{itemize}
\item Vertex operator of type I :$\Phi_{i}(v)$
\item Vertex operator of type II:$\Psi^{*}_{i}(v)$
\end{itemize}
\noindent These vertex operators realize the Faddeev-Zamolodchikov (ZF) 
algebra with dynamical R-matrix as its structure coefficients 
\begin{eqnarray}
& &\Phi_{\nu}(v_2)\Phi_{\mu}(v_1)=R^{\ \ \mu'\nu'}_{F\mu\nu}
(v_1-v_2,\hat{\pi})\Phi_{\mu'}(v_1)\Phi_{\nu'}(v_2)\\
& &\Psi^{*}_{\mu}(v_1)\Psi^{*}_{\nu}(v_2)=-R^{*\mu\nu}_{F\mu'\nu'}
(v_1-v_2,\hat{\pi})\Psi^{*}_{\nu'}(v_2)\Psi^{*}_{\mu'}(v_1)\\
& &\Phi_{\nu}(v_1)\Psi^{*}_{\mu}(v_2)=\tau(v_1-v_2)\Psi^{*}_{\mu}(v_2) 
\Phi_{\nu}(v_1)
\end{eqnarray}
\noindent Let us introduce the other basic operators 
\begin{eqnarray*}
& &\eta_{1}(v)=e^{-i\sqrt{\frac{r-1}{r}}(Q-i2vlnxP)}:e^{-\sum_{m\neq 0}
\frac{\beta_{m}}{m}x^{-2vm}}:\\
& &\eta'_{1}(v)=e^{i\sqrt{\frac{r}{r-1}}(Q-i2vlnxP)}:e^{\sum_{m\neq 0}
\frac{\beta'_{m}}{m}x^{-2vm}}:\\
& &\xi (v)=e^{i\sqrt{\frac{2(r-1)}{r}}(Q-i2vlnxP)}:e^{\sum_{m\neq 0}
\frac{x^{m}+x^{-m}}{m}\beta_{m}x^{-2vm}}:\\
& &\xi '(v)=e^{-i\sqrt{\frac{2r}{r-1}}(Q-i2vlnxP)}:e^{-\sum_{m\neq 0}
\frac{x^{m}+x^{-m}}{m}\beta'_{m}x^{-2vm}}:
\end{eqnarray*}
\noindent where the q-deformed bosonic oscillators $\beta_{m},P,Q$ are
defined in Eq.(3).Then, the bosonization of  vertex operators are given
by[17,30,33]
\begin{eqnarray}
& &\Phi_{+}(v)=\eta_{1}(v)\ \ \ ,\ \ \ 
\Phi_{-}(v)=\oint_{C} \frac{d(x^{2v_1})}{2\pi ix^{2v_1}}\eta_{1}(v)\xi (v_1)
f(v_1-v,\hat{\pi})\\
& &\Psi_{+}(v)=\eta'_{1}(v)\ \ \ ,\ \ \ 
\Psi_{-}(v)=\oint_{C'} \frac{d(x^{2v_1})}{2\pi ix^{2v_1}}\eta'_{1}(v)\xi' (v_1)
f'(v_1-v,\hat{\pi})
\end{eqnarray}
\noindent where the integration contour $C$ is a simple closed curves around the 
origin satisfying $|xx^{2v}|<|x^{2v_1}|<|x^{-1}x^{2v}|$;$C'$ is chosen in
such a way that the poles $x^{2v-1+2n(r-1)}$ $(0\leq n)$are inside and the poles
$x^{2v+1-2n(r-1)}$ $(0\leq n)$ are outside, and 
\begin{eqnarray*}
f(v,w)=\frac{[v+\frac{1}{2}-w]_{r}}{[v-\frac{1}{2}]_{r}}\ \ ,\ \ 
f'(v,w)=\frac{[v-\frac{1}{2}+w]_{r-1}}{[v+\frac{1}{2}]_{r-1}}
\end{eqnarray*}
\noindent Set 
\begin{eqnarray*}
\Psi^{*}_{\mu}(v)=\Psi_{-\mu}(v)\frac{1}{[\hat{\pi}]_{r-1}}
\end{eqnarray*}

\noindent From the normal order relations given in appendix A, one can 
check that the bosonic realization for $\Phi_{i}$ and $\Psi^{*}$ in 
Eq.(52)---Eq.(53) satisfy ZF algebra in Eq.(49)---Eq.(51)[30,33].
\subsection{The  realization of algebra 
$A_{q,p;\hat{\pi}}(\widehat{gl_2})$ at level one by Miki's construction}
Let us introduce Miki's construction[34]
\begin{eqnarray}
& &L^{+\mu}_{\ \ \nu}(v,\hat{\pi})=\Phi_{\mu}(v)\Psi^{*}_{\nu}(v-\frac{1}{2})\\
& &L^{-\mu}_{\ \ \nu}(v,\hat{\pi})=\Phi_{\mu}(v-\frac{1}{2})\Psi^{*}_{\nu}(v)
\end{eqnarray}
\noindent Using the relations of ZF algebra in Eq.(49)---Eq.(51), one can 
prove the $L^{\pm}$-operators constructed above satisfy the definition of  
the elliptic quantum algebra Eq.(24)---Eq.(26) with $c=1$. Moreover, we have

\noindent {\bf Proposition 4.} The $L^{\pm}$-operators and two type vertex 
operators satisfy the folllowing relations 
\begin{eqnarray}
& &R_{F}^{+}(v_1-v_2,\hat{\pi})L_{1}^{+}(v_1,\hat{\pi})\Phi_{2}(v_2)=
\Phi_{2}(v_2)L^{+}_{1}(v_1,\hat{\pi})\\
& &R_{F}^{-}(v_1-v_2-\frac{1}{2},\hat{\pi})L_{1}^{-}(v_1,\hat{\pi})\Phi_{2}(v_2)
=\Phi_{2}(v_2)L^{-}_{1}(v_1,\hat{\pi})\\
& &L^{+}_{1}(v_1,\hat{\pi})\Psi^{*}_{2}(v_2)=
\Psi^{*}_{2}(v_2)L^{+}_{1}(v_1,\hat{\pi})R_{F}^{*+}(v_1-v_2-\frac{1}{2},\hat{\pi})\\
& &L^{-}_{1}(v_1,\hat{\pi})\Psi^{*}_{2}(v_2)=
\Psi^{*}_{2}(v_2)L^{-}_{1}(v_1,\hat{\pi})R_{F}^{*+}(v_1-v_2,\hat{\pi})
\end{eqnarray}
\noindent The proof is direct by using Miki's construction of 
$L^{\pm}$-operators  and ZF algebra Eq.(49)---Eq.(51).

From the Proposition 4, one can see that the vertex operators of ABF model 
are the intertwing operators of the elliptic algebra  
$A_{q,p;\hat{\pi}}(\widehat{gl_2})$ at level one, which satisfy 
some generized relations of q-affine algebra and its intertwing operators[31].

\section{The scaling limit algebra $A_{\hbar,\eta;\hat{\pi}}(\widehat{gl_2})$}
Another deformed Virasoro algbera---$\hbar$-Virasoro algebra can be considered
as symmetries of the massive integrable field theories[14] ,and at the
semi-classical level corresponds to the center of the Yangian Double with
 center $DY_{\hbar}(\widehat{gl_2})$ at the critical level[18]. In another
 way,$\hbar$-Virasoro algbera can be considered as the scaling limit of the
 q-Virasoro algebra[14]
\begin{eqnarray*}
x^{2v}=p^{-\frac{i\beta}{\hbar}}\ \ \ ,\ \ \ q=p^{-\frac{1}{\eta}}\ \ \
p\longrightarrow 1
\end{eqnarray*}
\noindent Moreover, the screening currents of $\hbar$-Virasoro
 algebra satisfies a closed algebra relations, which also can be considered
 as the scaling limit of that of q-Virasoro algbera in Eq.(9)-Eq.(15)[14].
Therefore, we can construct a generalizing algbera
  $A_{\hbar,\eta;\hat{\pi}}(\widehat{gl_2})$  as the scaling limit of algebra
  $A_{q,p;\hat{\pi}}(\widehat{gl_2})$, which is expected to be the
  symmetric algebra of k-fused Restricted sine-Gordon model and its
  (Drinfeld) current algebra would be the algebra of screening currents
  of some $\hbar$-deformed extended Virasoro algbera.Similarly, the algebra
  $A_{\hbar,\eta;\hat{\pi}}(\widehat{gl_2})$ can be formated in the dynamical
  $RLL=LLR^{*}$ form with the dynamical R-matrix being the trigonometric
  solution to the dynamical Yang-Baxter equation[14].

In this section, we restrict ourselves to the trigonometric dynamical
R-amtrix (or the scaling limit of R-matrix in Eq.(16) and Eq.(20)).
Without confusion with 
that in former section, we choose the same symbols as that of section 3.1 
\begin{eqnarray}
& &R_{F}(v,\hat{\pi})\equiv R_{F}(v,\hat{\pi},\eta)=\left(
\begin{array}{llll}a&&&\\&b&c&\\&d&e&\\&&&a\end{array}
\right)
\end{eqnarray}
\noindent and
\begin{eqnarray}
& &a(\beta,\hat{\pi})=\kappa(\beta)=exp\{\int^{\infty}_{0}
\frac{2sh\frac{\hbar t}{2}sh\frac{\hbar t}{2\eta}shi\beta t}
{sh\hbar t sh\frac{(1+\eta)\hbar t}{2\eta}}\frac{dt}{t}\ \ \ \ \ \ 
\}\\
& &\frac{b(\beta,\hat{\pi})}{a(\beta,\hat{\pi})}
=\frac{sin\pi\eta(\frac{i\beta}{\hbar})sin\pi\eta(\hat{\pi}-1)}
{sin\pi\eta(\hat{\pi})sin\pi\eta(\frac{i\beta}{\hbar}+1)}\ \ ,\ \
\frac{c(\beta,\hat{\pi})}{a(\beta,\hat{\pi})}
=\frac{sin\pi\eta sin\pi\eta(\frac{i\beta}{\hbar}+\hat{\pi})}
{sin\pi\eta(\hat{\pi})sin\pi\eta(\frac{i\beta}{\hbar}+1)}\\
& &\frac{d(\beta,\hat{\pi})}{a(\beta,\hat{\pi})}
=\frac{sin\pi\eta sin\pi\eta(-\frac{i\beta}{\hbar}+\hat{\pi})}
{sin\pi\eta(\hat{\pi})sin\pi\eta(\frac{i\beta}{\hbar}+1)}\ \ ,\ \ 
\frac{e(\beta,\hat{\pi})}{a(\beta,\hat{\pi})}
=\frac{sin\pi\eta(\frac{i\beta}{\hbar})sin\pi\eta(\hat{\pi}+1)}
{sin\pi\eta(\hat{\pi})sin\pi\eta(\frac{i\beta}{\hbar}+1)}\\
& &R^{\pm}_{F}(\beta,\hat{\pi})=\tau^{\pm}(\beta)R_{F}(\beta,\hat{\pi})
\ \ ,\ \ \tau^{+}(\beta)=ctg(\frac{i\pi\beta}{2\hbar})\ \ ,\ \
\tau^{-}(\beta)=-tg(\frac{i\pi\beta}{2\hbar})\nonumber
\end{eqnarray}
Define
\begin{eqnarray*}
R_{F}^{*\pm}(\beta,\hat{\pi})\equiv R^{\pm}_{F}(\beta,-\hat{\pi})|_
{\eta\longrightarrow \eta'}\ \ \ ,\ \ \frac{1}{\eta'}=\frac{1}{\eta}-c
\end{eqnarray*}
\noindent The algebra $A_{\hbar,\eta;\hat{\pi}}(\widehat{gl_2})$ is
generated by the matrices elements of $L^{\pm}$-operator which satisfy
the following relations
\begin{eqnarray}
& &R_{F}^{+}(\beta_1-\beta_2-\frac{i\hbar c}{2},\hat{\pi})
L^{+}_{1}(\beta_1,\hat{\pi})L^{-}_{2}(\beta_2,\hat{\pi})=
L^{-}_{2}(\beta_2,\hat{\pi})L^{+}_{1}(\beta_1,\hat{\pi})
R_{F}^{*+}(\beta_1-\beta_2+\frac{i\hbar c}{2},\hat{\pi})\\
& &R_{F}^{-}(\beta_1-\beta_2+\frac{i\hbar c}{2},\hat{\pi})
L^{-}_{1}(\beta_1,\hat{\pi})L^{+}_{2}(\beta_2,\hat{\pi})=
L^{+}_{2}(\beta_2,\hat{\pi})L^{-}_{1}(\beta_1,\hat{\pi})
R_{F}^{*-}(\beta_1-\beta_2-\frac{i\hbar c}{2},\hat{\pi})\\
& &R_{F}^{\pm}(\beta_1-\beta_2,\hat{\pi})L^{\pm}_{1}(\beta_1,\hat{\pi})
L^{\pm}_{2}(\beta_2,\hat{\pi})=L^{\pm}_{2}(\beta_2,\hat{\pi})
L^{\pm}_{1}(\beta_1,\hat{\pi})R_{F}^{*\pm}(\beta_1-\beta_2,\hat{\pi})\\
& &
\hat{\pi} L^{(\pm)\mu}_{\ \ \ \nu}(\beta,\hat{\pi})=
L^{(\pm)\mu}_{\ \ \ \nu}(\beta,\hat{\pi})(\hat{\pi}+(\nu \frac{1}{\eta}-
\frac{1}{\eta'}\mu))
\end{eqnarray}
Let
\begin{eqnarray*}
L^{\pm}(\beta,\hat{\pi})=
\left(\begin{array}{cc}1&0\\E^{\pm}(\beta)&1\end{array}\right)
\left(\begin{array}{cc}K^{\pm}_{1}(\beta)&0\\&K^{\pm}_{2}(\beta)
\end{array}\right)
\left(\begin{array}{cc}1&F^{\pm}(\beta)\\0&1\end{array}\right)
\end{eqnarray*}
\noindent be the Gauss decomposition of $L^{\pm}$-operators. 
For the convenience, we also introduce the following symbols 
\begin{eqnarray*}
& &R_{F}^{\pm}(\beta,\hat{\pi})=\left(\begin{array}{llll}
a^{\pm}(\beta)&&&\\&b^{\pm}(\beta)&c^{\pm}(\beta)&\\&d^{\pm}(\beta)
&e^{\pm}(\beta)&\\&&&a^{\pm}(\beta)\end{array}
\right)\ \ ,\ \  R_{F}^{*\pm}(\beta,\hat{\pi})=\left(\begin{array}{llll}
a'^{\pm}(\beta)&&&\\&b'^{\pm}(\beta)&c'^{\pm}(\beta)&\\&d'^{\pm}(\beta)
&e'^{\pm}(\beta)&\\&&&a'^{\pm}(\beta)\end{array}
\right)  
\end{eqnarray*}
\noindent Define the total currents $E(\beta)$ and  
$F(\beta)$ by the corresponding Ding-Frenkel correspondence
\begin{eqnarray}
& &E(\beta)=E^{+}(\beta)-E^{-}(\beta-\frac{i\hbar c}{2})\ \ \ ,\ \ \ 
F(\beta)=F^{+}(\beta-\frac{i\hbar c}{2})-F^{-}(\beta)
\end{eqnarray}
\noindent Substituting the Gauss decomposition of $L^{\pm}$-operator, we can
obtain the similar commutation relations bewteen $E(\beta)$,$F(\beta)$ and
$K^{\pm}_{i}(\beta)$ as those in the proposition 1, where the matrices elements of
R-matrix are in Eq.(60)---Eq.(63).

The algebra $A_{\hbar,\eta}(\widehat{gl_2})$ as the scaling limit of
the elliptic algebra $A_{q,p}(\widehat{gl_2})$ was studied by Khoroshkin
et al through the method of Gauss decomposition[8].The commutation relations
 of $E(\beta)$,$F(\beta)$ and $K^{\pm}_{i}(\beta)$ were obtained, which are
 quite different from ours.This is due to that they chose a nondynamical
 $``RLL=LLR^{*}"$ formalism. Consequently, the commutation
 relations of theirs dos not depend on the dynamical variable.If we introduce
 $H^{\pm}(\beta)$ as Eq.(69) and Eq.(70) ,which quite differ from that of
 Khoroshkin et al, algebra $A_{\hbar,\eta;\hat{\pi}}(\widehat{gl_2})$ and
 algebra $A_{\hbar,\eta}(\widehat{gl_2})$   share the same subalgebra---
 the (Drinfeld) current algebra of each one generated by $E(\beta),F(\beta),
 H^{\pm}(\beta)$.Although they share the same subalgebraic commutation
 relations,they yet relate to the different $E,F,K^{\pm}_{i}$ (or the
 different algebra) and consequently are associated with different vertex
 operators[18,14].Moreover, the two algebra are related with the different
 models (face type model for the dynamical algebra and vertex model for the
 nondynamical algebra) especially for the rational $\eta$ case.

Set
\begin{eqnarray}
& &H^{+}(\beta)
=K^{-}_{2}(\beta-\frac{i\hbar c}{4})\frac{sin\pi\eta'\hat{\pi}sin\pi\eta'}
{\pi\eta'sin\pi\eta'(\hat{\pi}-1)}K^{-}_{1}(\beta-\frac{i\hbar c}{4})^{-1}\\
& &H^{-}(\beta)
=K^{+}_{2}(\beta-\frac{i\hbar c}{4})\frac{sin\pi\eta'\hat{\pi}sin\pi\eta'}
{\pi\eta'sin\pi\eta'(\hat{\pi}-1)}K^{+}_{1}(\beta-\frac{i\hbar c}{4})^{-1}
\end{eqnarray}
\noindent we have the current algebra of algebra
$A_{\hbar,\eta;\hat{\pi}}(\widehat{gl_2})$ generated by $E(\beta),F(\beta)$
and $H^{\pm}(\beta)$ which is the
scaling limit of those of the elliptic algebra
$A_{q,p;\hat{\pi}}(\widehat{gl_2})$ .

\noindent {\bf Proposition 5.} The current algebra of algebra  
$A_{\hbar,\eta;\hat{\pi}}(\widehat{gl_2})$ is generated by
$E(\beta),F(\beta),H^{\pm}(\beta)$ and satisfies the following relations
\begin{eqnarray}
& &E(\beta_1)E(\beta_2)=\frac
{sin\frac{i\pi\eta}{\hbar}(\beta_1-\beta_2-i\hbar)}
{sin\frac{i\pi\eta}{\hbar}(\beta_1-\beta_2+i\hbar)}E(\beta_2)E(\beta_1)\\
& &F(\beta_1)F(\beta_2)=\frac
{sin\frac{i\pi\eta'}{\hbar}(\beta_1-\beta_2+i\hbar)}
{sin\frac{i\pi\eta'}{\hbar}(\beta_1-\beta_2-i\hbar)}F(\beta_2)F(\beta_1)\\
& &[E(\beta_1),F(\beta_2)]=\hbar\{
\delta(\beta_1-\beta_2-\frac{i\hbar c}{2})H^{+}(\beta_1-\frac{i\hbar c}{4})
-\delta(\beta_1-\beta_2+\frac{i\hbar c}{2})H^{-}(\beta_1+\frac{i\hbar c}{4})\}\\
& &H^{\pm}(\beta_1)E(\beta_2)=\frac
{sin\frac{i\pi\eta}{\hbar}(\beta_1-\beta_2-i\hbar \mp\frac{i\hbar c}{4})}
{sin\frac{i\pi\eta}{\hbar}(\beta_1-\beta_2+i\hbar \mp\frac{i\hbar c}{4})}
E(\beta_2)H^{\pm}(\beta_1)\\
& &H^{\pm}(\beta_1)F(\beta_2)=\frac
{sin\frac{i\pi\eta'}{\hbar}(\beta_1-\beta_2+i\hbar \pm\frac{i\hbar c}{4})}
{sin\frac{i\pi\eta'}{\hbar}(\beta_1-\beta_2-i\hbar \pm\frac{i\hbar c}{4})}
F(\beta_2)H^{\pm}(\beta_1)\\
& &H^{\pm}(\beta_1)H^{\pm}(\beta_2)=\frac
{sin\frac{i\pi\eta'}{\hbar}(\beta_1-\beta_2+i\hbar)
sin\frac{i\pi\eta}{\hbar}(\beta_1-\beta_2-i\hbar)}
{sin\frac{i\pi\eta'}{\hbar}(\beta_1-\beta_2-i\hbar)
sin\frac{i\pi\eta}{\hbar}(\beta_1-\beta_2+i\hbar)}
H^{\pm}(\beta_2)H^{\pm}(\beta_1)\\
& &H^{+}(\beta_1)H^{-}(\beta_2)=\frac
{sin\frac{i\pi\eta'}{\hbar}(\beta_1-\beta_2+i\hbar-\frac{i\hbar c}{2})
sin\frac{i\pi\eta}{\hbar}(\beta_1-\beta_2-i\hbar+\frac{i\hbar c}{2})}
{sin\frac{i\pi\eta'}{\hbar}(\beta_1-\beta_2-i\hbar-\frac{i\hbar c}{2})
sin\frac{i\pi\eta}{\hbar}(\beta_1-\beta_2+i\hbar+\frac{i\hbar c}{2})}
H^{-}(\beta_2)H^{+}(\beta_1)\\
\end{eqnarray}
\noindent It can be seen that when $c=1$ (i.e at level one), the current
algebra of algebra $A_{\hbar,\eta;\hat{\pi}}(\widehat{gl_2})$
be the algebra of the screening currents for $\hbar$-Virasoro algebra[14].
For higer level,it would be the algebra of screening currents for
$\hbar$-deformed extented Virasoro algebra.
Moreover, there exist the following
relations between the algebra $A_{\hbar,\eta;\hat{\pi}}(\widehat{gl_2})$ and
$A_{\hbar,\eta}(\widehat{gl_2})$ [8],and between the algebra
$A_{q,p;\hat{\pi}}(\widehat{gl_2})$ and $A_{q,p}(\widehat{gl_2})$[10] for
the generic $r$ and $\eta$ case
\begin{eqnarray*}
\begin{array}{lll}
A_{q,p;\hat{\pi}}(\widehat{gl_2})&\stackrel{{\bf\rm\tiny scaling\ \ limit}}
{\longrightarrow}&A_{\hbar,\eta;\hat{\pi}}(\widehat{gl_2})\\
\ \ \ \updownarrow {\bf\rm\tiny twisted}&&\ \ \
\updownarrow {\bf\rm\tiny twisted}\\
A_{q,p}(\widehat{gl_2})&\stackrel{{\bf\rm\tiny scaling\ \ limit}}
{\longrightarrow}&A_{\hbar,\eta}(\widehat{gl_2})
\end{array}
\end{eqnarray*}
\section{Discussions}
In this paper, we propose an elliptic algebra 
$A_{q,p;\hat{\pi}}(\widehat{gl_2})$ based on a dynamical relations 
$RLL=LLR^{*}$, where the dynamical R-matrix is of the $A^{(1)}_1$ type 
face model.The corresponding (Drinfeld) current algebra is the current
algebra generalizing the screening currents for q-Viarsoro algebra and
is a dynamical twisted algebra of $A_{q,p}(\widehat{gl_2})$, which can be
considered as the results of the correspondence between the $A^{(1)}_1$ face
model and the eight-vertex model with the generic $r$. Moreover,
the q-primary fields of q-Virasoro $\Phi_{\mu}$ and
 $\Psi^{*}_{\mu}$ which are also called as the vertex operators for 
 $A^{(1)}_1$ face model , are the intertwing operators of the elliptic 
 algebra $A_{q,p;\hat{\pi}}(\widehat{gl_2})$ at level one.

It is very interesting to extend the present formulation $RLL=LLR^{*}$ 
to the case of $A^{(1)}_{n-1}$. The corresponding elliptic algebra is   
$A_{q,p;\hat{\pi}}(\widehat{gl_n})$.The corresponding (Drinfeld) current
algebra of algebra $A_{q,p;\hat{\pi}}(\widehat{gl_n})$ would be the current
algebra generalizing the screening currents for q-deformed $W_{n}$ algebra, 
which is generated by $E_{j}(v),F_{j}(v)$ and $H^{\pm}_{j}(v)$ (j=1,...n-1)  
with the following relations 
\begin{eqnarray*}
& & E_i(v_1)E_j(v_2)=(-1)^{A_{ij}}\frac{[v_1-v_2-\frac{A_{ij}}{2}]_{r}}
{[v_1-v_2+\frac{A_{ij}}{2}]_{r}}E_j(v_2)E_i(v_1)\\
& & F_i(v_1)F_j(v_2)=(-1)^{A_{ij}}\frac{[v_1-v_2+\frac{A_{ij}}{2}]_{r-c}}
{[v_1-v_2-\frac{A_{ij}}{2}]_{r-c}}F_j(v_2)F_i(v_1)\\
& &[E_j(v_1),F_j(v_2)]=\frac{1}{x-x^{-1}}\{
\delta(v_1-v_2+\frac{c}{2})H_j^{+}(v_1+\frac{c}{4})
-\delta(v_1-v_2-\frac{c}{2})H_j^{-}(v_1-\frac{c}{4})\}\\
& &E_j(v_1)F_{j+1}(v_2)=-F_{j+1}(v_2)E_j(v_1)\ \ ,\ \ 
F_j(v_1)E_{j+1}(v_2)=-E_{j+1}(v_2)F_j(v_1)\\ 
& &[E_j(v_1),F_l(v_2)]=0\ \ \ ,\ \ \ |j-l|>1\\
& &H^{\pm}_{i}(v_1)E_{j}(v_2)=\frac{[v_1-v_2-\frac{A_{ij}}{2}\mp\frac{c}{4}]
_{r}}{[v_1-v_2+\frac{A_{ij}}{2}\mp\frac{c}{4}]_{r}}E_j(v_2)H^{\pm}_j(v_1)\\
& &H^{\pm}_{i}(v_1)F_{j}(v_2)=\frac{[v_1-v_2+\frac{A_{ij}}{2}\pm\frac{c}{4}]
_{r-c}}{[v_1-v_2-\frac{A_{ij}}{2}\pm\frac{c}{4}]_{r-c}}F_j(v_2)
H^{\pm}_j(v_1)\\
& &H^{\pm}_i(v_1)H^{\pm}_j(v_2)=\frac{[v_1-v_2-\frac{A_{ij}}{2}]_{r}[v_1-v_2
+\frac{A_{ij}}{2}]_{r-c}}
{[v_1-v_2+\frac{A_{ij}}{2}]_{r}[v_1-v_2
-\frac{A_{ij}}{2}]_{r-c}}H^{\pm}_j(v_2)H^{\pm}_i(v_1)\\
& &H^{\pm}_i(v_1)H^{\mp}_j(v_2)=
\frac{[v_1-v_2-\frac{A_{ij}}{2}\mp\frac{c}{2}]_{r}[v_1-v_2
+\frac{A_{ij}}{2}\pm\frac{c}{2}]_{r-c}}
{[v_1-v_2+\frac{A_{ij}}{2}\mp\frac{c}{2}]_{r}[v_1-v_2
-\frac{A_{ij}}{2}\pm\frac{c}{2}]_{r-c}}H^{\mp}_j(v_2)H^{\pm}_i(v_1)
\end{eqnarray*}
\noindent where $H^{-}_{j}(v)=H^{+}_{j}(v+\frac{c}{2}-r)$ and the matrix
$A_{ij}$ is the Cartan matrix for $A^{(1)}_{n-1}$ Lie algebra 
\begin{eqnarray*}
A_{ij}=2\delta_{ij}-\delta_{i+1,j}-\delta_{i-1,j}
\end{eqnarray*}
\noindent The above algberaic relations could be derived by the Gauss 
decomposition of the $L^{\pm}$-operators corresponding to the 
dynamical R-matrix of $A^{(1)}_{n-1}$ face model. We will present the  
results in the further paper.

{\it \bf Acknowledgements.} Some results of this paper was given as a talk 
by one of authors B.Y.Hou in the XIIth international congress of
mathematical physics 13-19 July, 1997, Brisbane, Australia and in the
international workshop on ``Statistical Mechanics and Integrable Systems"
28 July-8 August 1997, Canberra,Australia.He would like to thank 
Prof.A.J.Bracken, M.D.Gould, R.J.Baxter and B.M.M.McCoy 
for their hospitality.He also would like to thank Prof.V.E.Korepin, H.Awata
and H.Konno for their interest about this problem in Australia.
One of authors W.L.Yang would like to thank
Prof.K.J.Shi, L.Zhao and Dr.H.Fan for their helpful discussions.

{\it \bf Note added:} After our paper was submitted to the electronic archive,
Prof. H.Konno informed us that " Recently, Jimbo has succeeded to derive
the algebra $U_{q,p}(\widehat{sl_2})$ by the Gauss decomposition of a central
extended dynamical RLL---relations with the R-matrix introduced by Enriquez
 and Felder. This results allows us to clarify a Hopf algebra structure of
 $U_{q,p}(\widehat{sl_2})$. This work along this line is now in progress".

\section*{Appendix}
\subsection*{A.The normal order relation for basic operators}
The normal order relations for the screening currents $E(v)$ and $F(v)$ of
q-deformed Virasoro algebra are
\begin{eqnarray*}
& &E(v_1)E(v_2)=x^{\frac{4(r-1)v_1}{r}}t(v_2-v_1):E(v_1)E(v_2):\\
& &F(v_1)F(v_2)=x^{\frac{4rv_1}{r-1}}t'(v_2-v_1):F(v_1)F(v_2):\\
& &E(v_1)F(v_2)=\frac{x^{-4v_1}}{(1-xx^{2(v_2-v_1)})(1-x^{-1}
x^{2(v_2-v_1)})}:E(v_1)F(v_2):\\
& &F(v_2)E(v_1)=\frac{x^{-4v_2}}{(1-xx^{2(v_1-v_2)})(1-x^{-1}
x^{2(v_1-v_2)})}:E(v_1)F(v_2):
\end{eqnarray*}
\noindent The normal order relations for the basic operators in type I
and type II vertex operators read as
\begin{eqnarray*}
& &\eta_{1}(v_1)\eta_{1}(v_2)=x^{\frac{(r-1)v_1}{r}}g_{1}(v_2-v_1)
:\eta_{1}(v_1)\eta_{1}(v_2):\\
& &\eta_{1}(v_1)\xi(v_2)=x^{\frac{2(1-r)v_1}{r}}s(v_2-v_1)
:\eta_{1}(v_1)\xi(v_2):\\
& &\xi(v_2)\eta_{1}(v_1)=x^{\frac{2(1-r)v_2}{r}}s(v_1-v_2)
:\eta_{1}(v_1)\xi(v_2):\\
& &\xi(v_1)\xi(v_2)=x^{\frac{4(r-1)v_1}{r}}t(v_2-v_1):\xi(v_1)\xi(v_2):\\
& &\eta'_{1}(v_1)\eta'_{1}(v_2)=x^{\frac{rv_1}{r-1}}g'_{1}(v_2-v_1)
:\eta'_{1}(v_1)\eta'_{1}(v_2):\\
& &\eta'_{1}(v_1)\xi'(v_2)=x^{\frac{-2rv_1}{r-1}}s'(v_2-v_1)
:\eta'_{1}(v_1)\xi'(v_2):\\
& &\xi'(v_2)\eta'_{1}(v_1)=x^{\frac{-2rv_2}{r-1}}s'(v_1-v_2)
:\eta'_{1}(v_1)\xi'(v_2):\\
& &\xi'(v_1)\xi'(v_2)=x^{\frac{4rv_1}{r-1}}t'(v_2-v_1):\xi'(v_1)\xi'(v_2):\\
& &\eta_{1}(v_1)\xi'(v_2)=(x^{2v_1}-x^{2v_2}):\eta_1(v_1)\xi'(v_2):\\
& &\xi'(v_2)\eta_{1}(v_1)=(x^{2v_2}-x^{2v_1}):\eta_1(v_1)\xi'(v_2):\\
& &\eta'_{1}(v_1)\xi(v_2)=(x^{2v_1}-x^{2v_2}):\eta'_1(v_1)\xi(v_2):\\
& &\xi(v_2)\eta'_{1}(v_1)=(x^{2v_2}-x^{2v_1}):\eta'_1(v_1)\xi(v_2):\\
& &\xi(v_1)\xi'(v_2)=\frac{x^{-4v_1}}{(1-xx^{2(v_2-v_1)})(1-x^{-1}
x^{2(v_2-v_1)})}:\xi(v_1)\xi'(v_2):\\
& &\xi'(v_2)\xi(v_1)=\frac{x^{-4v_2}}{(1-xx^{2(v_1-v_2)})(1-x^{-1}
x^{2(v_1-v_2)})}:\xi(v_1)\xi'(v_2):
\end{eqnarray*}
\noindent where
\begin{eqnarray*}
& &g_{1}(v)=\frac{\{x^{2+2v}\}\{x^{2+2r+2v}\}}{\{x^{4+2v}\}\{x^{2r+2v}\}}\ \ ,
\ \ s(v)=\frac{(x^{2r-1+2v};x^{2r})}{(x^{1+2v};x^{2r})}\\
& &t(v)=(1-x^{2v})\frac{(x^{2+2v};x^{2r})}{(x^{2r-2+2v};x^{2r})}\\
& &g'_{1}(v)=\frac{\{x^{2v}\}'\{x^{2+2r+2v}\}'}{\{x^{2r+2v}\}'
\{x^{2+2v}\}'}\ \ ,\ \ \{z\}'=(z;x^{2(r-1)},x^{4})\\
& &s'(v)=\frac{(x^{2r-1+2v};x^{2(r-1)})}{(x^{-1+2v};x^{2(r-1)})}\ \ ,\ \ 
t'(v)=(1-x^{2v})\frac{(x^{-2+2v};x^{2(r-1)})}{(x^{2r+2v};x^{2(r-1)})}
\end{eqnarray*}
\subsection*{B.The proof of the commutation relations between $K^{\pm}_{i}(v)$,
$E(v)$ and $F(v)$}

\noindent The proof is direct substitution the Gauss decomposition 
of $L^{\pm}$-operators in the relations Eq.(24)---Eq.(26)(It should be
careful to deal with the order between the dynamical R-matrices
and the Guass components of $L^{\pm}$-operators).Here, we give the proof
of Eq.(34) as an example. After some straightforward calculatation, one can
obtain the following commutation relations between the partial currents
$E^{\pm}(v)$ and $F^{\pm}(v)$
\begin{eqnarray*}
& &[E^{\pm}(v_2),F^{\pm}(v_1)]=
K^{\pm}_{1}(v_1)^{-1}\frac{c^{\pm}(v_1-v_2)}{b^{\pm}(v_1-v_2)}
K^{\pm}_{2}(v_1)-K^{\pm}_{2}(v_2)\frac{c'^{\pm}(v_1-v_2)}{b'^{\pm}(v_1-v_2)}
K^{\pm}_{1}(v_2)^{-1}\\
& &[E^{-}(v_2),F^{+}(v_1)]=
K^{+}_{1}(v_1)^{-1}\frac{c^{+}(v_1-v_2+\frac{c}{2})}{b^{+}(v_1-v_2+\frac{c}{2})}
K^{+}_{2}(v_1)-K^{-}_{2}(v_2)\frac{c'^{+}(v_1-v_2-\frac{c}{2})}
{b'^{+}(v_1-v_2-\frac{c}{2})}K^{-}_{1}(v_2)^{-1}\\
& &[E^{+}(v_2),F^{-}(v_1)]=K^{-}_{1}(v_1)^{-1}\frac{c^{-}(v_1-v_2-\frac{c}{2})}
{b^{-}(v_1-v_2-\frac{c}{2})}K^{-}_{2}(v_1)-K^{+}_{2}(v_2)
\frac{c'^{-}(v_1-v_2+\frac{c}{2})}{b'^{-}(v_1-v_2+\frac{c}{2})}
K^{+}_{1}(v_2)^{-1}
\end{eqnarray*}
\noindent Then, we have 
\begin{eqnarray*}
& &[E(v_1),F(v_2)]=
K^{-}_{2}(v_1+\frac{c}{2})\{
\frac{c'^{+}(v_1-v_2+\frac{c}{2})}{b'^{+}(v_1-v_2+\frac{c}{2})}
-\frac{c'^{-}(v_1-v_2+\frac{c}{2})}{b'^{-}(v_1-v_2+\frac{c}{2})}\}
K^{-}_{1}(v_1+\frac{c}{2})^{-1}\\
& &\ \ \ \ \ +K^{+}_{2}(v_1)\{
\frac{c'^{-}(v_1-v_2-\frac{c}{2})}{b'^{-}(v_1-v_2-\frac{c}{2})}
-\frac{c'^{+}(v_1-v_2-\frac{c}{2})}{b'^{+}(v_1-v_2-\frac{c}{2})}\}
K^{+}_{1}(v_1)^{-1}
\end{eqnarray*}
\noindent Using the following identity
\begin{eqnarray*}
& &
\frac{c'^{+}(v+i\epsilon)}{b'^{+}(v+i\epsilon)}
-\frac{c'^{-}(v-i\epsilon)}{b'^{-}(v-i\epsilon)}
=\frac{1}{x-x^{-1}}\delta(v)\frac{[\hat{\pi}]_{r-c}[1]_{r-c}}{\theta'_{r-c}
[\hat{\pi}-1]_{r-c}}\ \ \ {\rm When}\ \ \epsilon\longrightarrow 0\\
& &\ \ \ \theta'_{t}=(x-x^{-1})\frac{\partial}{\partial v}[v]_{t}|_{v=0}
\end{eqnarray*}
\noindent we can get the Eq.(34).
\subsection*{C. The proof of the commutation relations Eq.(37)---Eq.(43).}
First we prove the relations Eq.(37) and Eq.(38).From the properties Eq.(27) 
and Eq.(28), using the Gauss decomposition for $L^{\pm}$-operators we have 
\begin{eqnarray*}
& &R^{\pm}(v_1,\hat{\pi})E(v_2)=E(v_2)R^{\pm}(v_1,\hat{\pi}-2c)\ \ ,\ \ 
R^{\pm}(v_1,\hat{\pi})F(v_2)=F(v_2)R^{\pm}(v_1,\hat{\pi})\\ 
& &R^{*\pm}(v_1,\hat{\pi})E(v_2)=E(v_2)R^{*\pm}(v_1,\hat{\pi})\ \ ,\ \ 
R^{*\pm}(v_1,\hat{\pi})F(v_2)=F(v_2)R^{*\pm}(v_1,\hat{\pi}-2c)
\end{eqnarray*}
\noindent Notice that the relations Eq.(33a) and Eq.(33b) , we have 
\begin{eqnarray*}
& &E(v_1)E(v_2)=\frac{[v_1-v_2-1]_{r}}{[v_1-v_2+1]_{r}}E(v_2)E(v_1)\\
& &F(v_1)F(v_2)=\frac{[v_1-v_2+1]_{r-c}}{[v_1-v_2-1]_{r-c}}F(v_2)F(v_1)
\end{eqnarray*}
\noindent In order to obtain the relations between $H^{\pm}(v)$, the 
reltions between $H^{\pm}(v)$ and $E(v),F(v)$ , we should deal with 
Eq.(30b). It can be rewritten as  
\begin{eqnarray*}
\frac{a^{\pm}(v_1-v_2)}{a'^{\pm}(v_1-v_2)}
K^{\pm}_{2}(v_2)\frac{[v_1-v_2+1]_{r-c}[\hat{\pi}]_{r-c}}
{[v_1-v_2]_{r-c}[\hat{\pi}-1]_{r-c}}K^{\pm}_{1}(v_1)^{-1}=
K^{\pm}_{1}(v_1)^{-1}\frac{[v_1-v_2+1]_{r}[\hat{\pi}]_{r}}
{[v_1-v_2]_{r}[\hat{\pi}-1]_{r}}K^{\pm}_{2}(v_2)
\end{eqnarray*}
\noindent Taking the limit of $v_1\longrightarrow v_2$ in both side of the  
above equation, we have 
\begin{eqnarray*}
K^{\pm}_{2}(v)\frac{[1]_{r-c}[\hat{\pi}]_{r-c}}
{\theta'_{r-c}[\hat{\pi}-1]_{r-c}}K^{\pm}_{1}(v)^{-1}=
K^{\pm}_{1}(v)^{-1}\frac{[1]_{r}[\hat{\pi}]_{r}}
{\theta'_{r}[\hat{\pi}-1]_{r}}K^{\pm}_{2}(v)
\end{eqnarray*}
\noindent Then we have two equivalent definition of $H^{\pm}(v)$ 
\begin{eqnarray}
& &H^{+}(v)
=K^{-}_{2}(v+\frac{c}{4})\frac{[\hat{\pi}]_{r-c}[1]_{r-c}}
{\theta'_{r-c}[\hat{\pi}-1]_{r-c}}K^{-}_{1}(v+\frac{c}{4})^{-1}\nonumber\\
& &\ \ \ \ \ \ \ \ =K^{-}_{1}(v+\frac{c}{4})^{-1}
\frac{[\hat{\pi}]_{r}[1]_{r}}
{\theta'_{r}[\hat{\pi}-1]_{r}}K^{-}_{2}(v+\frac{c}{4})\\
& &H^{-}(v)
=K^{+}_{2}(v+\frac{c}{4})\frac{[\hat{\pi}]_{r-c}[1]_{r-c}}
{\theta'_{r-c}[\hat{\pi}-1]_{r-c}}K^{+}_{1}(v+\frac{c}{4})^{-1}\nonumber\\
& &\ \ \ \ \ \ \ \ =K^{+}_{1}(v+\frac{c}{4})^{-1}
\frac{[\hat{\pi}]_{r}[1]_{r}}{\theta'_{r}[\hat{\pi}-1]_{r}}
K^{+}_{2}(v+\frac{c}{4})
\end{eqnarray}
\noindent From these two equivalent definitions of $H^{\pm}(v)$ , one can 
check  Eq.(40)---Eq.(41). Here, we give the proof of Eq.(42) as an example
\begin{eqnarray*}
& &H^{+}(v_1)H^{+}(v_2)=
K^{-}_{2}(v_1+\frac{c}{4})
\frac{[\hat{\pi}]_{r-c}[1]_{r-c}}{\theta'_{r-c}[\hat{\pi}-1]_{r-c}}
K^{-}_{1}(v_1+\frac{c}{4})^{-1}K^{-}_{1}(v_2+\frac{c}{4})^{-1}
\frac{[\hat{\pi}]_{r}[1]_{r}}{\theta'_{r}[\hat{\pi}-1]_{r}}
K^{-}_{2}(v_2+\frac{c}{4})\\
& & \ \ =K^{-}_{2}(v_1+\frac{c}{4})
\frac{[\hat{\pi}]_{r-c}[1]_{r-c}}{\theta'_{r-c}[\hat{\pi}-1]_{r-c}}
\frac{a'^{-}(v_1-v_2)}{a^{-}(v_1-v_2)}
K^{-}_{1}(v_2+\frac{c}{4})^{-1}K^{-}_{1}(v_1+\frac{c}{4})^{-1}
\frac{[\hat{\pi}]_{r}[1]_{r}}{\theta'_{r}[\hat{\pi}-1]_{r}}
K^{-}_{2}(v_2+\frac{c}{4})\\
& & \ \ =\frac{[v_2-v_1]_{r-c}[1]_{r-c}}{\theta'_{r-c}[v_2-v_1+1]_{r-c}}
a^{-}(v_2-v_1)K^{-}_{2}(v_1+\frac{c}{4})\frac{1}{b'^{-}(v_2-v_1)}
K^{-}_{1}(v_2+\frac{c}{4})^{-1}\\
& &\ \ \ \ \ \ \ \ \ \ \ \ \ \times K^{-}_{1}(v_1+\frac{c}{4})^{-1}
\frac{[\hat{\pi}]_{r}[1]_{r}}{\theta'_{r}[\hat{\pi}-1]_{r}}
K^{-}_{2}(v_2+\frac{c}{4})\\
& & \ \ =
\frac{[v_2-v_1]_{r-c}[v_1-v_2]_{r}[1]_{r-c}[1]_{r}}
{\theta'_{r-c}\theta'_{r}[v_2-v_1+1]_{r-c}[v_1-v_2+1]_{r}}
K^{-}_{1}(v_2+\frac{c}{4})^{-1}\frac{a^{-}(v_2-v_1)}{b^{-}(v_2-v_1)}
K^{-}_{2}(v_2+\frac{c}{4})\\
& &\ \ \ \ \ \ \ \ \ \ \ \ \ \times K^{-}_{2}(v_1+\frac{c}{4})
\frac{a'^{-}(v_1-v_2)}{b'^{-}(v_1-v_2)}K^{-}_{1}(v_1+\frac{c}{4})^{-1}\\
& &\ \ =\frac{[v_1-v_2-1]_{r}[v_1-v_2+1]_{r-c}}
{[v_1-v_2+1]_{r}[v_1-v_2-1]_{r-c}}H^{+}(v_2)H^{+}(v_1)
\end{eqnarray*}
\noindent here we have used the identity
\begin{eqnarray*}
\frac{a'^{\pm}(v_1-v_2)}{a^{\pm}(v_1-v_2)}=
\frac{a^{\pm}(v_2-v_1)}{a'^{\pm}(v_2-v_1)}
\end{eqnarray*}
\noindent Similarily, we can prove the other relations among $H^{\pm}(v)$. 
The following identites are very useful for the proof  
\begin{eqnarray*}
a^{+}(v)a^{-}(-v)=1\ \ \  ,\ \ \ a'^{+}(v)a'^{-}(-v)=1
\end{eqnarray*}
\newpage
\section*{References}
\begin{enumerate}
\item R.Baxter, Exactly solved models in statistical mechanics (Academic 
Press, New York, 1992). 
\item C.N.Yang, {\it Phys. Rev. Lett.}{\bf 19} (1967) 1312. 

\item V.G.Drinfeld, {\it Sov. Mayh. Dokl.} {\bf 32} (1985) 254.

\item M.Jimbo, {\it Lett. Math. Phys.}{\bf 10} (1985) 63; 
{\it Lett. Math. Phys.}{\bf 11} (1986) 247; {\it Comm. Math. 
Phys.} {\bf 102} (1986) 537.

\item L.D.Faddeev, N.Reshetikhin and L.A.Takhtajan, {\it Algbera and 
Analysis} {\bf 1} (1989) 118.

\item N.Reshetikhin and M.A.Semenov-Tian-Shansky, {\it Lett. Math. Phys. } 
{\bf 19} (1990) 133.

\item J.Ding and I.B.Frenkel, {\it Comm. Math. Phys.} {\bf 155} (1993) 277.

\item S.Khoroshkin, D.Lebedev and S.Pakuliak, {\it Phys. Lett.} {\bf A 222}
(1996) 381; Elliptic
algebra $A_{q,p}(\widehat{sl_2})$ in the scaling limit, {\it q-alg/9702002}. 

\item K.Iohara and M.Kohno,{\it Lett. Math. Phys.} {\bf 37} (1996) 319.

\item O.Foda, K.Iohara, M.Jimbo, R.Kedem, T.Miwa and H.Yan,{\it Lett. Math. 
Phys.} {\bf 32} (1994) 259; Notes on highest weight modules of the 
elliptic algebra $A_{q,p}(\widehat{sl_2})$ ,{\it Prog. Theor. Phys.} 
{\bf Suppl. 118} (1995) 1. 

\item J.Shiraishi, H.Kubo, H.Awata and S.Odake, {\it Lett. Math. Phys.} {\bf 
38} (1996) 33.

\item H.Awata, H.Kubo, S.Odake and J.Shiraishi, {\it Comm. Math. Phys.} {\bf 
179} (1996) 401; Virasoro-type symmetries in Solvable models.

\item B.Feigin and E.Frenkel, {\it Comm. Math. Phys.} {\bf 178} (1996) 653.

\item B.Y.Hou and W.L.Yang, A $\hbar$-deformed Virasoro algbera as the 
hidden symmetry algebra of the Restricted sine-Gordon model, 
{\it hep-th/9612235}, ( accepted by {\it Comm. Theo. Phys.}) (1997); 
A $\hbar$-deformed $W_N$ algebra and its 
vertex operators, {\it hep-th/9701101, Jour. Phys.} {\bf A30} (1997),6131 
;The quantum $\hbar$-deofrmed $W_N$ algera and the algebra   
of its screening currents,(in preparation).

\item A.A.Belavin, A.M.Polyakov and A.B.Zamolodchikov, {\bf Nucl. Phys.}
 {\bf B241}(1984) 333.

\item S.Lukyanov, {\it Phys. Lett }{\bf B367}(1996) 121.

\item S.Lukyanov and Y.Pugai, {\it Nucl. Phys. }{\bf B473}(1996) 631; {\it 
Jour. Exp. Theor. Phys.} {\bf 82}(1996) 1021.

\item X.M.Ding, B.Y.Hou and L.Zhao, $\hbar$-(Yangian) deformation Miura Map
and Virasoro algebra, {\it q-alg/9701014};
The algebra $A_{\hbar,\eta}(\widehat{g})$
and Hopf family of algebras.

\item G.Andrews, R.Baxter and J.Forrester, {\it Jour. Stat. Phys.} 
{\bf 35}(1984) 193.

\item M.Jimbo, M.Lashkevich, T.Miwa and Y.Pugai, Lukyanov's screening 
operators for the deformed Virasoro algbera, {\it hep-th/9607177}. 

\item B.Feigin, M.Jimbo, T.Miwa, A.Odesskiiad and Y.Pugai, Algbera of 
screening operators for the deformed $W_n$ algebra .

\item B.Enriquez and G.Felder, Elliptic quantum group $E_{\tau,\eta}(sl_2)$ 
and Quasi-hopf algebras,{\it q-alg/9703018}.

\item J.Avan, O.Babelon and E.Billey, {\it Comm. Math. Phys.} {\bf 178}
(1996) 281.

\item O.Babelon, D.Bernard and E.Billey, {\it Phys. Lett.} {\bf B375}(1996) 
89.

\item G.Felder, {\it Comm. Math. Phys.}{\bf 176}(1996) 133; Elliptic quantum
groups,{\it hep-th/9412207}.

\item E.Date, M.Jimbo, A.Kuniba, T.Miwa and M.Okado, {\it Nucl. Phys.} 
{\bf B290 [FS20]}(1987) 231; {\it Adv. Stud. in Pure Math.} {\bf 16}(1988) 
17.

\item H.Konno, q-deformation of the coset conformal field theory and the
fusion RSOS model, Talks given in the XIIth International Congress of
Mathematical Physics,13-19 July,1997, Brisbane, Australia; An Elliptic
algebra $U_{q,p}(\widehat{sl_2})$ and the Fusion RSOS Model,
{\it q-alg/9709013}.

\item B.Y.Hou, K.J.Shi and Z.X.Yang, {\it Jour. Phys. }{\bf A 26}(1993) 4951.

\item M.Jimbo, T.Miwa and M.Okado, {\it Nucl. Phys.} {\bf B300[FS22]}(1988) 
74 .

\item Y.Asai, M.Jimbo, T.Miwa and Y.Pugai, {\it Jour. Phys. }{\bf A29}(1996) 
6595.

\item I.B.Frenkel and N.Reshetikhin, {\it Comm. Math. Phys. }{\bf 146}
(1992) 1.

\item M.Jimbo and T.Miwa, Algebraic analysis of solvable lattice models,
RIMS(1994), 981.

\item H.Fan, B.Y.Hou, K.J.Shi and W.L.Yang, the elliptic quantum algebra 
$A_{q,p}(\widehat{sl_n})$ and its bosonization at level one, 
{\it hep-th/9704024, IMPNWU-970401};L.D.Faddeev,{\it Lett. Math. Phys.}
{\bf 34} (1995),249.

\item K.Miki, {\it Phys. Lett.}{ \bf A186}(1994), 217.

\item T.Miwa and R.Weston, Boundary ABF models,{\it hep-th/961004}.

\item M.Jimbo, H.Konno and T.Miwa, Massless XXZ model and degeneration of the
elliptic algbera $A_{q,p}(\widehat{sl_2})$, {\it RIMS-1105}.

\item B.Y.Hou and W.L.Yang, $\hbar$(Yangian) deformed Virasoro algebras as
a dynamically twisted algebra $A_{\hbar,\eta}(\widehat{gl_2})$,Talk given
in The XIIth International Congress of Mathematical Physics,13-19 July, 1997,
Brisbane,Australia;Talk given in The International Workshop on ``Statistical
Mechanics and Integrable Systems",28 July-8 August, 1997,Canberra,Australia.

\item J.Bagger, D.Nemeshansky and S.Yaukielowicz, {\it Phys. Rev. Lett.} {\bf
60} (1988) 389.

\item D.Friedan, Z.Qiu and S.Shenker, {\it Phys. Lett.} {\bf B 151} (1985) 37.

\item D.Gepner and Z.Qiu {\it Nucl. Phys.} {\bf B285} (1987) 423; F.Ravanini,
{\it Mod. Phys. Lett.} {\bf A } (1988) 397.
\end{enumerate}

\end{document}